%% file: DGP_FVCA8.tex
\documentclass[12pt,a4paper,twoside,final,notitlepage, leqno]{article}
\usepackage[english]{babel}
\usepackage[T1]{fontenc}  
\usepackage{epsfig, graphicx, amssymb}
\usepackage{amsmath,amsthm,epsfig,amsfonts}  
\usepackage{float}
\usepackage{color}
\usepackage{mathptmx}       
\usepackage{type1cm}        
\usepackage{makeidx}         
\usepackage{graphicx}        
\usepackage{amsmath, amsfonts, amssymb}
\usepackage{color}
\usepackage{tikz}

\newcommand{\moneq}{\vspace*{-7pt} \begin{equation} \displaystyle } 
\newcommand{\moneqstar}{\vspace*{-6pt} \begin{equation*} \displaystyle } 
\newcommand{\monendstar}{\vspace*{-6pt} \end{equation*}   }
\newcommand{\monend}{\vspace*{-7pt} \end{equation}   }

\newcommand{\dd}{\,{\rm d}}

\def\section*#1{}

\usepackage{fancyhdr}
\renewcommand{\headrulewidth}{0pt}

\usepackage{amscd}  
\usepackage{cleveref} 
\let\div\undefined
\DeclareMathOperator{\div}{div} 
\DeclareMathOperator{\Span}{Span}

\DeclareMathOperator{\supp}{Supp} 
\DeclareMathOperator{\cotan}{cotan} 
\DeclareMathOperator{\gT}{\nabla_{\T}} 
\newcommand{\R}{\mathbb R}
\newcommand{\Hdiv}{{\rm H}(\div,\Omega)}
\newcommand{\ld}  {{\rm L}^2(\Omega)} 
\newcommand{\huo} {{\rm H}^1_0(\Omega)} 
\newcommand{\ldd} {\left [ \ld \right]^2} 
\newcommand{\vphi}{\varphi}
\newcommand{\vphia}{\varphi_a}
\newcommand{\vphias}{\varphi_a^\star}
\newcommand  {\un}  {1\hspace{-0.09cm} {\rm l}}
\newcommand  {\ai}  {a\in \T^1_i, ~ \partial^c a = (K,L)}
\newcommand  {\ab}  {a\in \T^1_b, ~ \partial^c a =(K)}
 \newcommand  {\T}   {\mathcal{T}}
\newcommand  {\drtb}{ \left( \vphias \right)_{a\in \T^1}}

\newcommand  {\ca}{ ( \vphias,\vphia)_0}

\newcommand  {\di}    {\displaystyle}

\newcommand  {\nref}{\labelcref}

\begin{document} 

\fancypagestyle{plain}{ \fancyfoot{} \renewcommand{\footrulewidth}{0pt}}
\fancypagestyle{plain}{ \fancyhead{} \renewcommand{\headrulewidth}{0pt}}

~

  \vskip 2.1 cm

\centerline {\bf \LARGE  Raviart Thomas  Petrov-Galerkin Finite Elements  }

 \bigskip  \bigskip \bigskip

\centerline { \large   Fran\c{c}ois Dubois$^{ab}$,  Isabelle Greff$^{c}$ and Charles Pierre$^{c}$}

\smallskip  \bigskip 

\centerline { \it  \small   
$^a$  Dpt. of Mathematics, University Paris-Sud,  B\^at. 425, F-91405  Orsay, France.} 

\centerline { \it  \small   
$^b$    Conservatoire National des Arts et M\'etiers, LMSSC laboratory,  F-75003 Paris, France.} 

\centerline { \it  \small  $^c$   
Laboratoire de Math\'ematiques et de leurs Applications, 
CNRS, Universit\'e de Pau, France. }

\centerline {francois.dubois@u-psud.fr, isabelle.greff@univ-pau.fr, charles.pierre@univ-pau.fr}


\bigskip  

\centerline {  25 February 2017  
{\footnote {\rm  \small $\,$ Contribution presented at the 8th 
Finite Volumes for Complex Applications conference, Lille, 12-16 June 2017.
Published in {\it Springer Proceedings in Mathematics \& Statistics}, 
volume~199, ``Finite Volumes for Complex Applications VIII~-~Methods and Theoretical Aspects'' 
(Eds C.~Canc\`es and P.~Omnes), pages~341-349, 2017, ISBN: 978-3-319-57396-0,  
DOI 10.1007/978-3-319-57397-7\_27. 
Edition 12~December~2017. }}}



\bigskip 
 \bigskip  \bigskip  
\noindent {\bf \large Abstract} 

\noindent 
The general theory of  Babu\v{s}ka 
  ensures  necessary and sufficient conditions for a
  mixed problem in classical or Petrov-Galerkin form
  to be  well posed in the sense of Hadamard.
  Moreover, the mixed method of Raviart-Thomas 
  with low-level elements can be interpreted as a finite volume method with a
non-local gradient. 
In this contribution, we propose a variant of type
  Petrov-Galerkin to ensure a
  local computation of the gradient at the interfaces of the elements.
  The in-depth study of stability leads to a specific choice
  of the test functions. With this choice, we show on the one hand that the
  mixed Petrov-Galerkin obtained is identical to the finite volumes scheme 
``volumes finis \`a 4 points'' (``VF4'')
 of Faille, Gallo\"uet  and   Herbin 
  and to the  condensation of mass approach  developed
  by Baranger, Maitre and Oudin. 
On the other hand, we show  the stability via an inf-sup condition 
and finally the convergence with the usual methods of mixed finite elements.

\bigskip \bigskip   \newpage \noindent {\bf \large    1) \quad  Introduction  }    

\fancyhead[EC]{\sc{Fran\c cois  Dubois, Isabelle Greff and Charles Pierre}} 
\fancyhead[OC]{\sc{Raviart Thomas  Petrov-Galerkin Finite Elements}} 
\fancyfoot[C]{\oldstylenums{\thepage}}

\smallskip \noindent 
{\bf Discrete gradient}

\smallskip  \noindent 
In the sequel,  $\Omega\subset \R^2$ denotes an open bounded convex 
 with a polygonal boundary.
The functional spaces $\ld$, $\huo$ and $\Hdiv$ are considered. 
The  $L^2$-scalar products on $\ld$ and on $\ldd$ are similarly denoted  $(\cdot,\cdot)_0$, 
without ambiguity.
Being set a triangulation $\T$ of $\Omega$, $P^0$ and $RT$ denote the associated finite 
element spaces  of the piecewise constant functions on the mesh and  the 
Raviart Thomas vector fields of order 0  \cite{RT77}, precise definitions follow 
in    Section~2. 
\\
The two unbounded operators,
gradient 
\moneqstar 
  \nabla:~  \ld \supset  \huo \rightarrow \ldd 
\monendstar 
 and divergence   
\moneqstar  
 \div:~  \ldd  \supset \Hdiv  \rightarrow  \ld 
\monendstar 
together satisfy the Green formula:  for $u\in\huo$ and $p\in\Hdiv$:
\moneqstar 
  (\nabla  u, p)_0 = -(u, \div p)_0
\monendstar 
Identifying $\ld$ and $\ldd$ with their topological dual spaces using the $L^2$-scalar 
product yields the following property, 
\moneqstar 
\nabla  = -\div ^\star \,, 
\monendstar 
that is a weak definition of the gradient on $\huo$.

\smallskip \noindent 
We search to define a \textit{discrete gradient} denoted $\gT$ on $P^0$ also based on a 
similar weak formalism. 
Starting from the divergence operator 
\moneqstar 
  \div:~  RT \rightarrow P^0 \,,
\monendstar 
one can define 
$ \,\,   \div^\star:~  \left(P^o\right)' \rightarrow \left(RT\right)' , 
\,\, $
between the algebraic dual spaces of $P^0$ and $RT$ respectively.
The natural basis for $P^0$ is made of the indicator functions of the mesh triangles,
that is orthogonal for the $L^2$-scalar product.
Therefore, $P^0$ is identified with its algebraic dual space $\left(P^o\right)'$.   
On the contrary, the Raviart Thomas basis $\left\{\vphi_a, \,  a \in \T^1\right\}$ of $RT$ 
(denoting by $\T^1$ the mesh edge set,   see Section~2) 
has no orthogonality property 
and cannot be used directly (see below)
to identify $RT$ with $\left(RT\right)'$.             
For this reason, a general identification process of $\left(RT\right)'$  
to a subspace $RT^\star\subset \Hdiv$ so that,
\moneqstar  
  RT^\star = \Span\left( \vphi_a^\star, \quad  a\in\T^1 \right) \,,
\monendstar 
with,
\moneq
  \label{eq:RT*-0}
  \vphi_a^\star \in \Hdiv, \quad 
  (\vphi_a^\star,\vphi_a)_0 \neq 0,
\monend
and the orthogonality property,
\moneq
  \label{eq:RT*}
  (\vphi_a^\star,\vphi_b)_0 = 0 \quad {\rm for} \quad a,b\in\T^1, \quad a\neq b,
\monend
is considered.
Setting $\Pi:~ RT\rightarrow  RT^\star$ with  $\Pi \vphi_a = \vphi_a^\star$, 
we have the following diagram, and general definition for the discrete gradient,
\moneq
\label{eq:grad-discret}
\begin{CD} 
RT         & @>{\div}>>    & P_0    
 \\
@ V{\Pi}VV & \hfill        & @VV{id}V
\\
RT^\star   & @ <<{\div^\star}< & P_0
\end{CD}
\quad, \qquad \gT = -\Pi^{-1}\div^\star~~:~P^0\rightarrow RT.
\monend
The definition of the discrete gradient is effective 
once $\ \, \{  \vphi_a^\star, \, a\in\T^1 \, \} \, $ has been set.
Various choices are possible.
The first choice is to set $RT^\star=RT$, and therefore to build 
$\left\{ \vphi_a^\star, \,  a\in\T^1 \right\}$   
with a Gram Schmidt orthogonalization process on the Raviart Thomas basis.
Such a choice has an important drawback. The  dual base function $\vphi_a^\star$ 
does not conserve a support located around the edge $a$. The discrete gradient 
matrix will be a full matrix related with the Raviart Thomas mass matrix inverse. 
This is not relevant with regard to the original gradient operator that is local in space.
This choice corresponds to the classical mixed finite element discrete gradient 
that is known to be associated with a full matrix.
In order to overcome this problem, Baranger, Maitre and Oudin \cite{BMO96} proposed to 
lump the mass matrix of the mixed finite element method. 
By doing this, they obtain a discrete {\it local} gradient. 
\\
A second choice, proposed in Thomas-Trujillo \cite{TT1} and also by 
one of us in \cite{D00, D02, Du10},  
that will be investigated in this paper, is to search for a dual basis satisfying, 
in addition to the orthogonality property \nref{eq:RT*}, the localization constraint,
  \moneq
    \label{eq:RT*-2}
    \forall ~a\in\T^1,\quad 
    \supp(\vphi^\star_a) \subset  \supp (\vphi_a),
  \monend
  in order to impose locality to the discrete gradient. 
With such a constraint 
the discrete gradient of $u\in P^0$  will be defined on each edge $a\in \T^1$
 only from the two values of $u$ on each side of $a$.
  In this context it is no longer asked to have $\vphi^\star_a \in RT$ so that 
$RT \neq RT^\star$: 
thus, this is a Petrov-Galerkin discrete formalism, 
as defined and used
{\it a priori} in the article of  Babu{\v{s}}ka \cite {Ba71}.   

 \bigskip \bigskip   \noindent {\bf \large    2) \quad  
Background and notations}
\label{section:notations}

\smallskip \noindent 
{\bf Meshes} 


\smallskip \noindent 
A conformal triangle mesh $\T$ of $\Omega$ in the sense followed by 
is considered.
The angle, vertex, edge and triangle sets of $\T$ are respectively denoted 
$\T^{-1}$, $\T^{0}$, $\T^1$ and $\T^2$.  
For $K\in\T^2$ (\textit{resp.} $a\in\T^1$) its area (\textit{resp.} length) 
is denoted $\vert  K\vert $ (\textit{resp.}  $\vert  a\vert $).

\noindent 
Let $K\in\T^2$.
Its three edges are denoted $a_{K,i}$, the unit normal to $a_{K,i}$ 
pointing outwards $K$ is denoted $n_{K,i}$, 
Its three vertices and angles are denoted $W_{K,i}$ and  $\theta_{K,i}$ respectively, 
so that $W_{K,i}$ and  $\theta_{K,i}$ are opposite to $a_{K,i}$ (see Fig.~1).

\noindent 
Let $a\in\T^1$.
One of its two unit normal is chosen and denoted $n_a$. This sets  an orientation for $a$.
Let $S_a,~ N_a$ be the two vertices of $a$, ordered so that $(n_a, S_aN_a)$ has a direct orientation.
\\
The sets $\T^1_i$ and $\T^1_b$ of the internal and boundary edges respectively are defined as,
\moneqstar  
  \T^1_b = \left\{ a\in\T^1, \quad  a\subset \partial\Omega\right\}, \,\, 
  \T^1_i = \T^1 - \T^1_b \, . 
\monendstar 
Let $a\in\T^1_i$. Its coboundary $\partial^c a$ is made of a unique ordered pair $K$, 
$L\in\T^2$ so that $a\subset  \partial K \cap \partial L$ and so that $n_a$ points 
from $K$ towards $L$. In such a case the following notation will be used:
\moneqstar    
  \ai
\monendstar 
and we will denote $W_a$ (\textit{resp.} $E_a$) the vertex of K (\textit{resp.} $L$) opposite to $a$.
Let $a\in\T^1_b$, $n_a$ is assumed to point towards the outside of $\Omega$. 
Its coboundary is made of a single $K\in\T^2$ so that $a\subset \partial K$, 
which situation is denoted as follows:
\moneqstar   
  \ab
\monendstar 
and we will denote $W_a$ the vertex of K opposite to $a$.
If $a\in\T^1$ is an edge of $K\in \T^2$, the angle of $K$ opposite to $a$ is denoted $\theta_{a,K}$ 
 (see Fig.~2).
%

 \bigskip \centerline {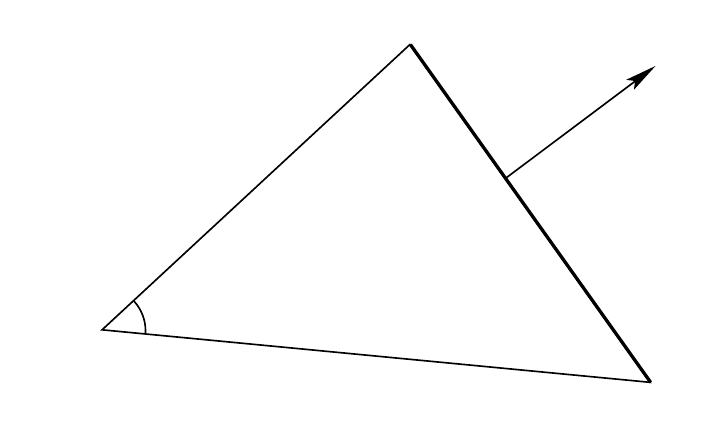}

\noindent
{\bf Fig. 1} Mesh notations for a triangle $K\in\T^2$  
  \label{fig:cell}

 \smallskip \centerline {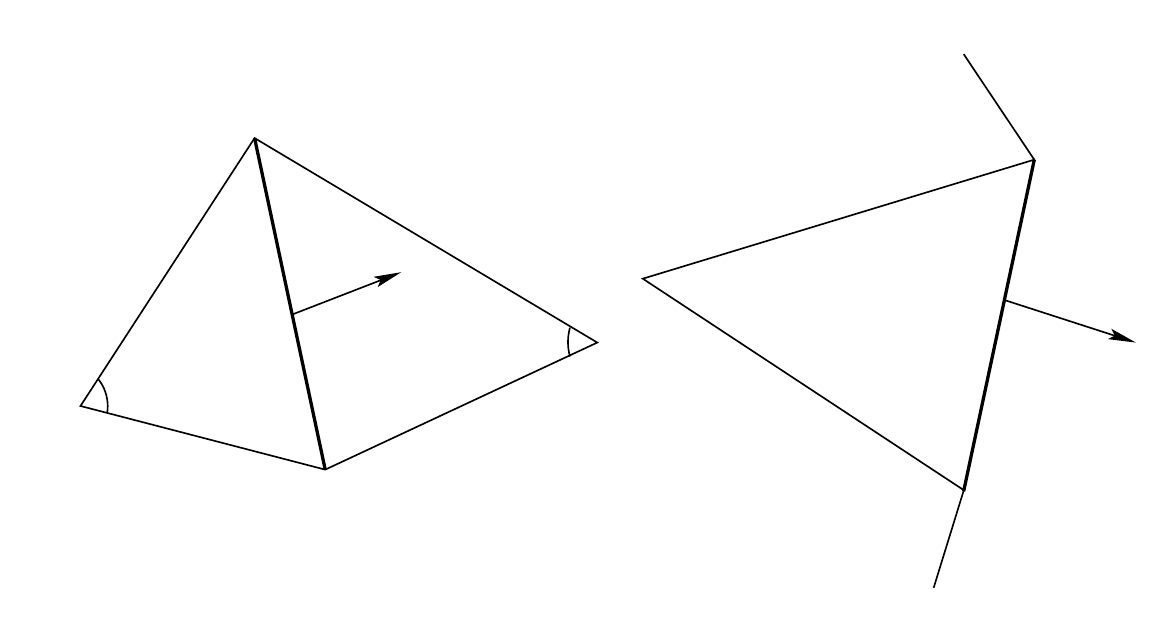}

\noindent
{\bf Fig. 2} Mesh notations for an internal edge (left) and for a boundary edge (right)
  \label{fig:edge}

\bigskip \bigskip \noindent 
{\bf The finite element spaces}  
\label{sec:fem}

\noindent 
Relatively to a mesh $\T$, the finite element spaces $P^0$ and $RT$  will be considered.
The space 
$P^0 \subset  \ld$ is the space of piecewise constant functions on the mesh triangles.
The indicators $\un_K$ for $K\in\T^2$ form a basis of $P^0$.
To $u\in P^0$ is associated the vector $(u_K)_{K\in\T^2}$ so that
\moneqstar 
  u = \sum_{K\in \T^2} u_K \, \un_K  \,. 
\monendstar 
The space 
$RT\subset  \Hdiv$ is the Raviart Thomas of order 0 finite element space introduced in~\cite{RT77}. 
An element $p\in RT$ is uniquely determined by its fluxes 
\moneqstar  
 \, p_a := \int_a p\cdot n_a \, \dd s \,\,\, 
 {\rm for} \,\,\,  \, a\in\T^1 \, .  
\monendstar 
The classical basis $\{\vphi_a, \, a \in \T^1\}$    
of $RT$ is so that 
\moneqstar 
 \, \int_b \vphi_a\cdot n_b \dd s = \delta_{ab} \,  
 \,\, \,\,  {\rm for} \,\, {\rm   all}  \, \, b\in\T^1  
\monendstar 
and 
with $\delta_{ab}$ the Kronecker symbol.
For each  $ \, p\in RT \, $ we can associate a discrete vector of  fluxes  
$(p_a)_{a\in\T^1}$ and we have 
$ \, p = \sum_{a\in\T^1} p_a \vphi_a $.

\noindent 
The \textit{local Raviart Thomas basis functions} are defined, 
for $K\in\T^2$ and $i=1$, 2, 3, by:
   \begin{equation}
    \vphi_{K,i}(x) = \dfrac{1}{2\vert  K\vert } \nabla \vert x-W_{K,i}\vert ^2
    \quad {\rm on} \quad  K 
   \quad  {\rm and} \quad  \vphi_{K,i} = 0 \quad {\rm otherwise}.
  \end{equation}
With that definition:
$ \,   \vphi_a = \vphi_{K,i} - \vphi_{L,j} \,\,$
if  $\, \ai \,\,  $
and 
$ \,\,  a=a_{K,i}=a_{L,j} \,\,$ 
$   \vphi_a = \vphi_{K,i} \, $ if $ \,\,  \ab\,\, $ 
and $ \,\, a=a_{K,i} \, $ 
and it is retrieved that $\supp(\vphi_a) =  K \cup L$ if $\ai$ or $\supp(\vphi_a) =  K$ in case $\ab$.
This provides a  second way to decompose $p\in RT$ as,
\moneqstar    
   p = \di\sum_{K\in\T^2} \sum_{i=1}^3 p_{K,i} \,\, \vphi_{K,i} \,,  
\monendstar 
where $p_{K,i} = \varepsilon p_a$ if $a=a_{K,i}$ with $\varepsilon = n_a \cdot n_{K,i} = \pm 1$.
Since $ \,\, \div \vphi_{K,i}=\frac{1}{|K|}$, 
the divergence operator $\div:~ RT \rightarrow P^0$ is given by,
\moneq
  \label{eq:RT-div}
  \div p = \sum_{K\in\T^2} \left(\div p\right)_K \un_K,
  \qquad  
  \left(\div p\right)_K= \dfrac{1}{| K|} \sum_{1=1}^3 p_{K,i}.
\monend

 \bigskip \bigskip   \noindent {\bf \large    3) \quad  
Raviart-Thomas dual basis }
\label{section:dualscheme}

\smallskip \noindent {\bf Definition 1} 
  \label{def:RT}

\noindent 
The family 
  $\drtb$ is said to be a Raviart Thomas dual basis if it satisfies (\ref{eq:RT*-0}), 
the orthogonality condition \nref{eq:RT*}, the localization condition \nref{eq:RT*-2} 
and the following \textit{flux normalization} condition:
\moneq
  \label{eq:RT*-3}
  \forall ~a,~b\in\T^1, \quad 
  \int_b \vphi_a^\star \cdot n_b  \dd s = \delta_{ab}, 
\monend
as for the Raviart Thomas basis functions $\vphia$, see  Section~2. 
%
In such a case, $RT^\star = \Span(\vphi_a^\star, \, a\in\T^1)$  
is the associated Raviart Thomas dual space, 
$\Pi:~\vphi_a\in RT \rightarrow  \vphi_a^\star \in RT^\star$ the projection onto $RT$  
and $\gT=-\Pi^{-1}\div^\star:~P^0\rightarrow RT$ the associated discrete gradient, 
as described in diagram \eqref{eq:grad-discret}.

\noindent 
The following algebraic relations will be useful.
From \cref{eq:RT*} one can check that, 
\moneq
  \label{eq:RT*.1}
  \forall ~p_1,~p_2\in RT,\quad 
  (\Pi p_1,p_2)_0 = ( p_1,\Pi p_2)_0 \, . 
\monend
The condition \eqref{eq:RT*-3} implies with the divergence theorem that, 
 $ \,   \forall ~p\in RT,\quad \forall ~K\in\T^2 $, 
$  \int_K \div p  \dd x = \int_K \div( \Pi p)  \dd x \,, $ 
and so that, 
\moneq  
  \label{eq:RT*-3.2}
  \forall ~(u,p)\in P^0\times  RT,\quad 
  (\div p, u)_0 = (\div (\Pi p), u)_0.
\monend
Now consider $u\in P^0$ and  $q\in RT^\star$. We have with (\ref{eq:RT*-3.2}),
$ \,   (u, \div q)_0 = (u, \div (\Pi^{-1} q))_0 = (\div^\star u, \Pi^{-1} q)_0 \, . $
Then with (\ref{eq:RT*.1}), 
$(u, \div q)_0 =(\Pi^{-1} (\div^\star u),  q)_0$.
As a result: 
\moneq
  \label{eq:RT*-1-3.1}
   \forall ~u\in P^0,\quad \forall ~q\in RT^\star,\quad 
  (u, \div q)_0 = -( \gT u, q)_0.
\monend

  \bigskip \noindent {\bf Proposition 1}  [Computation of the discrete gradient] 
\label{prop:disc-grad-prop-1}

\noindent 
Let $\drtb$ be a Raviart Thomas dual basis.
The discrete gradient is given 
for $u\in P^0$, by the
relation $\gT u = \di\sum_{a\in \T^1} p_a \vphi_a $ with the conditions
%
\moneq 
  \label{eq:disc-grad-1}
\left\{  \begin{array} {rcl}
\displaystyle 
\text{if} \,\, \ai, \, \,  p_a &=& \dfrac{u_L- u_K}{(\vphi_a,\vphi_a^\star)_0} \\ 
\displaystyle   \text{if} 
 \, \ab, \,\, p_a &=& \dfrac{- u_K}{(\vphi_a,\vphi_a^\star)_0} \, .
\end{array} \right.  \monend

\noindent This proposition deserves comments.
\\
The result of the localization condition (\ref{eq:RT*-2}) is, as expected, 
a local discrete gradient: its value on an edge $a\in \T^1$ only depends 
on the values of the scalar function $u$ on each sides of $a$.
The definition of the discrete gradient on the external edges implicitly 
takes into account a zero value for the scalar data $u$ on the domain boundary. 
This is relevant since the divergence with domain the full space $\Hdiv$ 
has for adjoint the gradient with domain $\huo$, which adjoint property 
has been translated at a discrete level. 
The formulation of the discrete gradient in proposition~1 
brings to the fore the coefficients $(\vphi_a^\star,\vphi_a)_0$: more details 
follow in the next subsection.
%

\bigskip \noindent
\textbf{Petrov-Galerkin discretization for the Dirichlet Poisson  problem}

\noindent 
Consider the following Dirichlet Poisson problem on $\Omega$,
\moneqstar   
 \,    - \Delta u = f \in \ld \,,  
 \,\,  u=0 \,\, {\rm  on } \, \,  \partial\Omega \, . 
\monendstar 
Consider a mesh $\T$ and  a Raviart Thomas dual basis $\drtb$.
Let us denote $V= P^0 \times RT$ and $V^\star= P^0 \times RT^\star$.
The  mixed
Petrov-Galerkin discretization of the Poisson problem 
is: find $(u,p)\in V$ so that,
\moneq
  \label{eq:PG-laplace}
  \forall ~ (v,q)\in V^\star,\quad 
  (p,q)_0 + (u,\div q)_0 = 0 
  \quad  \text{and} \quad 
  -(\div p, v)_0         = (f,v)_0  \, .
\monend
%
%
The mixed  Petrov-Galerkin discrete problem \nref{eq:PG-laplace} reformulates as:
find $(u,p)\in V$ so that,
\moneqstar  
\,   \forall ~ (v,q) \in V^\star, \,\,   
  B\bigl((u,p),(v,q)\bigl)=(f,v)_0 \,  
\monendstar 
where the bilinear  form $B$ 
is defined for 
$(u,p)\in V$ and  $ (v,q)\in  V^\star$ by,
\moneqstar   
   B\bigl((u,p),(v,q)\bigl) = (u,\div q)_0+(p,q)_0 - (\div p,v)_0  \, . 
\monendstar 
%
%
%

  \bigskip \noindent {\bf Proposition 2}   [Solution of the mixed discrete problem] 
  \label{prop:PG-disc}
 
\noindent 
The pair $(u,p)\in V$ is a  solution of problem \nref{eq:PG-laplace} if and only if 
\moneq
  \label{eq:FV-laplace}
  \gT u = p, \qquad  -\div( \gT u) = f_\T,
\monend
where $f_{\T}\in P^0$ is the projection of $f$, defined by,
\moneqstar 
   f_\T = \sum_{K\in \T^2} f_K \, \un_K,\,\,   
\,   f_K = \frac{1}{|K|} \int_K f \dd x \, .  
\monendstar 
If $\,\, (\vphi_a, \vphi_a^\star)>0 \,\, $ for all $ \,\, a\in\T^1$,  
then problem \nref{eq:PG-laplace} has a  unique solution.
%

\noindent 
Proposition~1 shows an equivalence between the 
mixed Petrov-Galerkin discrete problem  \nref{eq:PG-laplace} 
and the discrete problem (\ref{eq:FV-laplace}). Problem  (\ref{eq:FV-laplace}) 
actually is  a \textit{finite volume like} problem.
Precisely, 
it becomes: 
find $u\in P^0$ so that, for all $K\in \T^2$:
\moneqstar 
\di   \sum_{   \substack{
    a\in\T^1_i,~\partial^c a = (K,L)
    \\ ~~\text{or}~~ \partial^c a = (L,K)  
  }}
\dfrac{u_L - u_K}{\ca} \,\,  +
\sum_{
    a\in\T^1_b,~\partial^c a = (K)
  }
\frac{ - u_K}{\ca} = |K| \, f_K \, . 
\monendstar 
This \textit{finite volume like} problem only involves the coefficient $\ca$. 
We compute this scalar product  in the next section.
%
%

 \bigskip \bigskip   \noindent {\bf \large    4) \quad Retrieving the  ``VF4'' scheme}
%
%
%

\smallskip \noindent 
Let $g:~(0,1)\rightarrow \R$ so that,
\moneqstar 
 \,   \int_0^1 g  \dd s  = 1 \,,  
 \,    \int_ 0^1 g(s)s^2  \dd s = 0 \, \, {\rm  and } \, 
 \, g(s) = g(1-s) \, .  
\monendstar 
On a mesh $\T$ are defined $g_{K,i}:~ a_{K,i} \rightarrow \R$ for $K\in\T^2$ and $i=1$, 2, 3 as,
\moneqstar 
   g_{K,i}(x) = {{g(s)}\over{|a_{K,i}|}}  \quad  {\rm  for } \,\,\,  
\,   x \, = \, s \, S_{K,i} + (1-s) \, N_{K,i} \, . 
\monendstar  
For $K\in \T^2$ is denoted $\delta_K~:K \rightarrow \R$ a function that satisfies 
\moneqstar 
  \int_K \delta_K \dd x = 1  \,\, {\rm  and } \,  
 \,  \int_K \delta_K(x) \,\, |x - W_{K,i}|^2 \, \dd x = 0 \, 
 \,  {\rm  for } \,   i = 1 , ~2,~ 3. 
\monendstar  
To a family $(\psi_{K,i})$ of functions on $\Omega$ for $K\in\T^2$ and for $i=1$, 2, 3 
is associated the family $(\psi_a)_{a\in\T^1}$ so that,
\moneq 
  \label{eq:loc-basis}
\left\{  \begin{array} {rcl}
\displaystyle 
\text{if} \,\,  \ai \quad \text{and} \,\,\, a=a_{K,i}=a_{L,j} \,, \,\, 
  \psi_a &=& \psi_{K,i} - \psi_{L,j} \\ 
\displaystyle  \text{if} \,\,   \ab \,   \text{and}  \,\,\,  a=a_{K,i}  \,, \,\,   
 \psi_a &=& \psi_{K,i} \, . 
\end{array} \right.  \monend
%

\bigskip \noindent {\bf Theorem 1} [Error estimations] 
  \label{thm:VF4}

\noindent 
  Assume that the mesh $\T$ angles all satisfy $0< \theta_{K,i} < \pi/2$.
  Consider a family $(\vphi^\star_{K,i})$ of vector fields on $\Omega$ 
for $K\in\T^2$ and for $i=1$, 2, 3 that satisfy, independently on $i$, on $K$, 
\moneq 
    \label{eq:phi-star-def-1.1}
    \div \vphi^\star_{K,i} = \delta_K\,, \quad  
\vphi^\star_{K,i} = 0 \,\,\,  \text{ otherwise} 
\monend 
and, on $\partial K$,
\moneq 
    \label{eq:phi-star-def-1.2}
    \vphi^\star_{K,i}\cdot n = g_{K,i} \,\,  \text{on} \,\,  a_{K,i} \,, \quad 
 \vphi^\star_{K,i}\cdot n = 0  \,\,\,   \text{otherwise} \, .  
\monend 
%

\noindent 
  Let  $\drtb$ be constructed with \cref{eq:loc-basis}.
  Then $\drtb$ is a Raviart Thomas dual basis.
   The coefficients $\ca$ only depend on mesh $\T$ geometry, as follows
\moneq 
    \label{eq:coef-cotan}
\left\{  \begin{array} {rcl}
\displaystyle 
\text{ for } \ai \,\, \text{ then } \,\,  \ca &=& \left(\cotan \theta_{a,K}
      +\cotan \theta_{a,L}\right) / 2   \\ 
\displaystyle  \text{ for }  \ab \,\,  \text{ then } \,\, 
 \ca &=& \cotan \theta_{a,K} /2 \, . 
\end{array} \right.  \monend

\noindent 
  The mixed Petrov-Galerkin discrete problem  (\ref{eq:FV-laplace}) 
for the Poisson  equation 
has a unique solution and 
coincides with the classical ``VF4'' scheme introduced in \cite{FGH91}
(see also Faille \cite{Fa92} and Eymard {\it et al.} \cite{gallouet}).
%
%

\bigskip  \noindent 
Theorem~1 has various consequences.
Conditions in definition~1 
that must be satisfied by Raviart Thomas 
dual basis are replaced by sufficient conditions on $\delta_K$ and $g$.
In the sequel we will focus on such Raviart Thomas dual basis, though 
more general ones may exist: this will not be discussed in this contribution.
Assuming the existence of $g$ and $\delta_K$, 
the construction of such dual basis is very delicate. 
No explicit representation can {\it a priori} be obtained. 
Nevertheless,  a Raviart Thomas dual basis can be 
mathematically constructed by the following process. Consider 
$\vphi_{K,i}= \nabla u_{K,i}$ where $u_{K,i}$ is a solution of,
$ \,   \Delta u_{K,i} = \delta_K\quad \text{on}\quad K, \quad 
  \nabla  u_{K,i}\cdot n =  g_{K,i} \,\, $ on $ \,  a_{K,i} \,$
and 
$ \,  u_{K,i}\cdot n = 0 \, $ elsewhere on $ \, \partial K $. 
The compatibility condition for this problem is satisfied with 
the first statements 
and therefore $\vphi_{K,i}$  is well defined.

\noindent 
Whatever are $\delta_K$ and $g$, 
the coefficients $\ca$ will be unchanged: they only depend on the mesh geometry 
and are given by \cref{eq:coef-cotan}. 
Practically, this means that neither the $\drtb$ nor $\delta_K$ and $g$ need to be computed. 
The numerical scheme will always coincide with the ``VF4'' one.
Eventually, this provides a new point of view 
for the understanding and analysis of finite volume methods.
%
%

%



 \newpage 
 \bigskip \bigskip   \noindent {\bf \large    5) \quad Stability and convergence}
\label{sec:stab-conv}

\smallskip \noindent 
\textbf{General assumptions.}

\smallskip \noindent 
A couple of constant $0 < \theta_\star < \theta^\star <\pi /2$ is fixed and 
$\T$ will denote a mesh satisfying the uniform angle condition,
\moneq
  \label{eq:angle-cond}
  \forall ~ K\in\T^2,\quad i=1,~2,~3 :\quad 
  \theta_\star \le \theta_{K,i} \le \theta^\star.
\monend
%
Theorem 1 implies that the mixed Petrov-Galerkin discrete problem 
(\ref{eq:PG-laplace}) is independent on the particular choice made for  the Raviart Thomas dual basis. 
%

  \bigskip \noindent {\bf Theorem 2} [Error estimations]
\label{thm:conv-stab}

\noindent
There exists a constant $C$ independent on $\T$ and of $f$ 
in the Poisson problem 
so that the solution $(u_\T, p_\T)$ 
of the mixed  Petrov-Galerkin discrete problem (\ref{eq:PG-laplace}) satisfies,
\moneqstar  
   \Vert u_{\T} \Vert_0+\Vert p_{\T}\Vert_{\Hdiv}\leq  C \, \Vert f\Vert_{0} \, . 
\monendstar  
Denoting by $u$ the exact solution to the
Poisson problem 
and by $p=\nabla  u$ the following error estimates moreover holds,
\begin{equation}
  \label{eq:final-error}
 \Vert u-u_{\T} \Vert_0+\Vert p-p_{\T}\Vert_{\Hdiv}\leq 
  Ch_{\T}\Vert f\Vert_{1} \,, 
 \end{equation}
with $h_\T$ the mesh size.
%
%
%

\begin{proof}
  We first prove that the mixed Petrov-Galerkin formulation
  has a unique solution depending continuously on the data thanks to
  Babu\v{s}ka's work \cite{Ba71}.
The bilinear form $B$ 
is continuous on $V$: 
\moneqstar 
   | B(\xi,\eta)  |  \,\,\leq\,\, M \,\, \Vert \xi \Vert_{V}\, 
\Vert \eta \Vert_{V}\,,  \,  \forall \, \xi,\eta\in V \, . 
\monendstar 
The inf-sup stability condition 
relies   on a stability result \cite {D00, Du10, DGP-17}
and introduces a constant $ \, \beta > 0 \,$ such that for any mesh $ \,  {\T} $, 
\moneqstar 
\forall \, \xi \, \in P^0\times RT_0  \, {\rm such \, that} 
\,\,  \parallel \xi \parallel_{V} =  1 \,,  
 \,  \exists \, \eta \, \in  \, P^0\times RT_0^\star \,,
\, \,  \parallel \eta \parallel_{V}  \, \leq \, 1  \,  {\rm and }
 \, B(\xi,\,\eta) \, \geq \, \beta  \, .  
\monendstar 
The discrete ``infinity condition'' is satisfied \cite{D00}: 
\moneqstar 
 \,  \forall \eta \!  \in \! V^\star  \! \setminus \! \{0\},\!  
\,\,  \sup_{\xi \in V}  B(\xi,\eta)  \! = \! +\infty \, .   
\monendstar 
Then due to Babu\v{s}ka theorem valid also for Petrov-Galerkin
mixed formulation
the discrete scheme \eqref{eq:PG-laplace} has a unique solution
and 
\moneqstar 
 \, \Vert \xi - \xi_{_{\cal T}} \Vert_{V} \,\,\leq \,\, \big( 1+{M\over{\beta}} \big)
\inf_{\zeta\in V_{_{\cal T}}}\Vert \xi - \zeta \Vert_{V} \,. 
\monendstar   
In our case, this formulation is equivalent to
\moneq
  \label{eq:err-est}
\Vert u - u_{_{\cal T}}  \Vert _{0} \,+\,  \Vert p - p_{_{\cal T}}  \Vert_{\div}  
\,\leq\,  C \,\Big( \inf_{ v \in P^0} \, 
\Vert u-v\Vert_0 +\inf_{q\in RT}\Vert p-q \Vert_{\div} \Big)
\monend
for a constant $C=1+\frac{M}{\beta}$ dependent of $\T$ only through the lowest
and the highest angles $\theta_*$ and $\theta^*$.
We now precise an upper bound of the right-hand side of \eqref{eq:err-est}.
With the interpolation operators $\Pi_0:L^2(\Omega)\rightarrow P^0$
and $\Pi_{RT}:H^1(\Omega)^2\rightarrow RT^0$, we have 
\moneqstar 
 \, \Vert \, u - u_{_{\cal T}}  \, \Vert_{0} \,+\,  \Vert \, p - p_{_{\cal T}}  \, \Vert_{\div}
\,\leq\,  C \, \big(  \Vert \, u - \Pi_0   u \, \Vert_{0} 
+ \Vert \, p - \Pi_{\text{\tiny RT}}   p \, \Vert_{\div}  \big) \, . 
\monendstar 
On the other hand the interpolation errors are established by
Raviart and Thomas \cite{RT77} for the operator $\Pi_{\text{\tiny RT}} $:
%
\moneqstar  
 \,  \Vert \, u - \Pi_0   u \, \Vert_{0} 
 \, \leq \, h_{\cal T} \, \Vert u \, \Vert_{1}  \,,   
\monendstar 
\moneqstar  
\, \Vert \, p - \Pi_{\text{\tiny RT}}   p
\, \Vert_{0} \, \leq \, h_{\cal T} \, \Vert p \, \Vert _{1}  \,,   
\quad
\,   \Vert \, \div \big( p - \Pi_{\text{\tiny RT}}   p \big) 
  \, \Vert_{0} \, \leq \, h_{\cal T} \, \Vert \div p \, \Vert _{1} \, . 
\monendstar 
Then 
\moneqstar 
 \,  \Vert \, u - u_{_{\cal T}}  \, \Vert _{0} \,+\,  \Vert \, p - p_{_{\cal T}}  \, \Vert_{\div}  
\,\leq\,    C \,  h_{\cal T} \,  \big( \,  \Vert u \, \Vert  _{1}
+ \Vert p \, \Vert  _{1} + \Vert \div p \, \Vert  _{1} \,  \big) \, . 
\monendstar 

\smallskip \smallskip  \noindent 
Since
$ - \Delta u = f \,\,\, {\rm in} \, \Omega $,
with $f\in L^2(\Omega)$ and $\Omega$ convex, then
$u\in H^2(\Omega)$ and 
$\Vert u\Vert_2 \leq c \Vert f\Vert_0$.
Moreover $p=\nabla u$ and $\div p =-f$ leads to 
\moneqstar 
 \,  \Vert \, u - u_{_{\cal T}}  \, \Vert _{0} \,+\,  \Vert \, p - p_{_{\cal T}}  \, \Vert_{\div}  
\,\leq\,  C \,  h_{\cal T} \,  \big( 2 \Vert f \, \Vert_0
+ \Vert f \, \Vert  _{1}   \big) \, . 
\monendstar 
Finally, we get  
\moneqstar 
\,  \Vert \, u - u_{_{\T}}  \, \Vert _{0} \,+\,  \Vert \, p - p_{_{\cal T}}  \, \Vert_{\div}  
\,\leq\,  C \,  h_{\T}  \, \Vert  \, f \, \Vert_1   \,, 
\monendstar 
 that is exactly \eqref{eq:final-error}.
\end{proof}

 \bigskip  \bigskip    \noindent {\bf \large    6) \quad Possible extensions}
\label{sec:extensions}

\smallskip \noindent 
Our analysis for the Laplace equation is also {\it a priori} valid for three space dimensions. 
Moreover, the extension of the scheme to equations with tensorial coefficients 
is also possible in principle. 
To build a dual Raviart-Thomas basis for these problems
 is one of our objectives for a future contribution.


\bigskip \bigskip      \noindent {\bf  \large  References }   

\bibliographystyle{abbrv}
\bibliography{biblio}

\end{document}

%% file: cell_K.pdf_t
\begin{picture}(0,0)%
\includegraphics{cell_K.pdf}%
\end{picture}%
\setlength{\unitlength}{4144sp}%
\begingroup\makeatletter\ifx\SetFigFont\undefined%
\gdef\SetFigFont#1#2#3#4#5{%
  \reset@font\fontsize{#1}{#2pt}%
  \fontfamily{#3}\fontseries{#4}\fontshape{#5}%
  \selectfont}%
\fi\endgroup%
\begin{picture}(3226,1967)(278,-1228)
\put(1818,-419){\makebox(0,0)[lb]{\smash{{\SetFigFont{10}{12.0}{\familydefault}{\mddefault}{\updefault}{\color[rgb]{0,0,0}{\large $K$}}%
}}}}
\put(1021,-734){\makebox(0,0)[lb]{\smash{{\SetFigFont{10}{12.0}{\familydefault}{\mddefault}{\updefault}{\color[rgb]{0,0,0}$\theta_{K,i}$}%
}}}}
\put(3331,434){\makebox(0,0)[lb]{\smash{{\SetFigFont{10}{12.0}{\familydefault}{\mddefault}{\updefault}{\color[rgb]{0,0,0}$n_{K,i}$}%
}}}}
\put(451,-1006){\makebox(0,0)[lb]{\smash{{\SetFigFont{10}{12.0}{\familydefault}{\mddefault}{\updefault}{\color[rgb]{0,0,0}$W_{K,i}$}%
}}}}
\put(2386,-376){\makebox(0,0)[lb]{\smash{{\SetFigFont{10}{12.0}{\familydefault}{\mddefault}{\updefault}{\color[rgb]{0,0,0}$a_{K,i}$}%
}}}}
\end{picture}%

%% file: edge_a.pdf_t
\begin{picture}(0,0)%
\includegraphics{edge_a.pdf}%
\end{picture}%
\setlength{\unitlength}{4144sp}%
\begingroup\makeatletter\ifx\SetFigFont\undefined%
\gdef\SetFigFont#1#2#3#4#5{%
  \reset@font\fontsize{#1}{#2pt}%
  \fontfamily{#3}\fontseries{#4}\fontshape{#5}%
  \selectfont}%
\fi\endgroup%
\begin{picture}(5363,2853)(1988,-1896)
\put(2504,-869){\makebox(0,0)[lb]{\smash{{\SetFigFont{10}{12.0}{\familydefault}{\mddefault}{\updefault}{\color[rgb]{0,0,0}$\theta_{a,K}$}%
}}}}
\put(3839,-327){\makebox(0,0)[lb]{\smash{{\SetFigFont{10}{12.0}{\familydefault}{\mddefault}{\updefault}{\color[rgb]{0,0,0}$n_a$}%
}}}}
\put(2866,-667){\makebox(0,0)[lb]{\smash{{\SetFigFont{10}{12.0}{\familydefault}{\mddefault}{\updefault}{\color[rgb]{0,0,0}$K$}%
}}}}
\put(3800,-685){\makebox(0,0)[lb]{\smash{{\SetFigFont{10}{12.0}{\familydefault}{\mddefault}{\updefault}{\color[rgb]{0,0,0}$L$}%
}}}}
\put(3101,-186){\makebox(0,0)[lb]{\smash{{\SetFigFont{10}{12.0}{\familydefault}{\mddefault}{\updefault}{\color[rgb]{0,0,0}{\large $a$ }}%
}}}}
\put(6285,-799){\makebox(0,0)[lb]{\smash{{\SetFigFont{10}{12.0}{\familydefault}{\mddefault}{\updefault}{\color[rgb]{0,0,0}{\large $a$}}%
}}}}
\put(7032,-412){\makebox(0,0)[lb]{\smash{{\SetFigFont{10}{12.0}{\familydefault}{\mddefault}{\updefault}{\color[rgb]{0,0,0}$n_a$}%
}}}}
\put(5664,-476){\makebox(0,0)[lb]{\smash{{\SetFigFont{10}{12.0}{\familydefault}{\mddefault}{\updefault}{\color[rgb]{0,0,0}$K$}%
}}}}
\put(3101,394){\makebox(0,0)[lb]{\smash{{\SetFigFont{10}{12.0}{\familydefault}{\mddefault}{\updefault}{\color[rgb]{0,0,0}$N_a$}%
}}}}
\put(3499,-1412){\makebox(0,0)[lb]{\smash{{\SetFigFont{10}{12.0}{\familydefault}{\mddefault}{\updefault}{\color[rgb]{0,0,0}$S_a$}%
}}}}
\put(4843,-186){\makebox(0,0)[lb]{\smash{{\SetFigFont{10}{12.0}{\familydefault}{\mddefault}{\updefault}{\color[rgb]{0,0,0}$W_a$}%
}}}}
\put(4743,-799){\makebox(0,0)[lb]{\smash{{\SetFigFont{10}{12.0}{\familydefault}{\mddefault}{\updefault}{\color[rgb]{0,0,0}$E_a$}%
}}}}
\put(6460,652){\makebox(0,0)[lb]{\smash{{\SetFigFont{10}{12.0}{\familydefault}{\mddefault}{\updefault}{\color[rgb]{0,0,0}$\partial \Omega$}%
}}}}
\put(2161,-1636){\makebox(0,0)[lb]{\smash{{\SetFigFont{10}{12.0}{\familydefault}{\mddefault}{\updefault}{\color[rgb]{0,0,0} $\ai$}%
}}}}
\put(2206,-1051){\makebox(0,0)[lb]{\smash{{\SetFigFont{10}{12.0}{\familydefault}{\mddefault}{\updefault}{\color[rgb]{0,0,0}$W_a$}%
}}}}
\put(4321,-646){\makebox(0,0)[lb]{\smash{{\SetFigFont{10}{12.0}{\familydefault}{\mddefault}{\updefault}{\color[rgb]{0,0,0}$\theta_{a,L}$}%
}}}}
\put(4861,-1636){\makebox(0,0)[lb]{\smash{{\SetFigFont{10}{12.0}{\familydefault}{\mddefault}{\updefault}{\color[rgb]{0,0,0} $\ab$}%
}}}}
\end{picture}%